\documentclass[reqno]{amsart}
\usepackage[reset,a4paper,vmargin=5truecm,hmargin=4truecm]{geometry}
\usepackage{amsmath,amsthm,amssymb}

\usepackage{graphicx}
\usepackage{url}
\usepackage{cite}

\begin{document}

\title[Exact value of Tammes problem for $N=10$]{Exact value of 
Tammes problem for $N=10$}

\address{The Interdisciplinary Institute of Science, Technology and Art\\
Suzukidaini-building 211, 2-5-28 Kitahara, Asaka-shi, Saitama, 351-0036, 
Japan}
\email{ismsugi@gmail.com}

\address{The Institute of Statistical Mathematics, 10-3 Midori-cho, 
Tachikawa, Tokyo 190-8562, Japan}
\email{tanemura@ism.ac.jp}

\author{Teruhisa Sugimoto \and Masaharu Tanemura}
\maketitle	

\begin{abstract}
Let $C_{i}$ ($\,i=1,\ldots ,N\,$) be the $i$-th open spherical cap of angular 
radius $r$ and let $M_{i}$ be its center under the condition that none of the 
spherical caps contains the center of another one in its interior. We consider 
the upper bound, $r_{N} $, (not the lower bound !) of $r$ of the case in which 
the whole spherical surface of a unit sphere is completely covered with $N$ 
congruent open spherical caps under the condition, sequentially for 
$i=2,\ldots ,N-1\,$, that $M_{i}$ is set on the perimeter of $C_{i-1}$, and 
that each area of the set $(\cup _{\nu =1}^{i-1}C_{\nu })\cap C_{i}$ becomes 
maximum. In this paper, for $N = 10$, we found out that the solutions of 
our sequential covering and the solutions of the Tammes problem were strictly 
correspondent. Especially, we succeeded to obtain the exact value $r_{10}$ for 
$N = 10$.
\end{abstract}

\section{Introduction}

The circle on the surface of a sphere is called a spherical cap. Among the 
problems of packing on the spherical surface, the closest packing of congruent 
spherical caps is the most famous, and is usually known as the Tammes 
problem~\cite{Tammes_1930} . 
The details of Tammes problem are as follows: ``How must $N$ congruent 
non-overlapping spherical caps be packed on the surface of a unit sphere so 
that the angular diameter of spherical caps will be as great as possible?'' 
The Tammes problem is mathematically proved solutions were known for 
$N=1,\,\ldots , 12,$ and $24$\footnote{In recent years, Oleg Musin and 
Alexey Tarasov proved the solution of the Tammes problem for $N =13$ 
and $14$~\cite{Musin_Tarasov_2012, Musin_Tarasov_2015}.}
\cite{Danzer_1963, Fejes_1972, Teshi_Ogawa_2000}. 
The systematic method of attaining these solutions has not been given.
The exact values in the cases for $N=1,\,\ldots , 9, 11,12$ and $24$ 
are known, but only in the case for $N= 10$, the value is approximate 
range $[1.154479,1.154480]$ by Danzer~\cite{Danzer_1963}.
 
We considered the packing problem by a systematic method which is different 
from the approach by Danzer. As a result, we obtained the exact value of 
angular diameter of spherical caps in the Tammes problem for $N=10$ as 
following~\cite{Sugi_Tane_2007}.

\begin{equation}
\begin{split}
\label{eq_r10}
r_{10} & =
\tan ^{ - 1}\left( {\left( {\frac{4}{\sqrt 3 }\cos 
\left( {\frac{1}{3}\tan ^{ - 1}\left( {\frac{\sqrt 3 \sqrt {229} }{9}} 
\right)} \right) + 3} \right)^{\frac{1}{2}}} \right) \\ 
& \approx 1.1544798334192707378319618404230\cdots \quad \text{rad}\;%
\text{.}
\end{split}
\end{equation}

\section{Outline of Method}

Let us show our systematic method.
Suppose we have $N$ congruent open spherical caps with angular radius $r$ on 
the surface $\mathcal{S}$ of the unit sphere and suppose that these spherical caps 
cover the whole spherical surface without any gap under the condition that 
none of them contains the center of another one in its interior.
Further we suppose that the spherical caps are put on $\mathcal{S}$ sequentially 
in the manner which is described just below. Let $C_{i}$ be the $i$-th open 
spherical cap and let $M_{i}$ be its center ($i=1,\ldots ,N$).
Our problem is to calculate the upper bound of $r$ for the sequential
covering, such that $N$ congruent open spherical caps cover the whole
spherical surface $\mathcal{S}$ under the condition that $M_{i}$ is set
on the perimeter of $C_{i-1}$, and that each area of set $(\cup _{\nu
=1}^{i-1}C_{\nu })\cap C_{i}$ becomes maximum in sequence for $i=2,\ldots ,N-1$. 
Here, we define a \textit{half-cap} as the spherical cap whose angular radius 
is $\frac{r}{2}$ and which is concentric with that of the original cap.
Let us suppose the centers of $N$ half-caps are placed on the positions 
of the centers of spherical caps $C_i$ ($i=1,\ldots ,N$) which are 
considered. At this time, we get the packing with $N$ congruent half-caps. 
Therefore, our sequential covering is in connection with the packing 
problem~\cite{Sugi_Tane_2006, Sugi_Tane_2007}.

\renewcommand{\figurename}{{\small Figure.}} 
\begin{figure}[tbp]
\centering\includegraphics[width=12cm,clip]{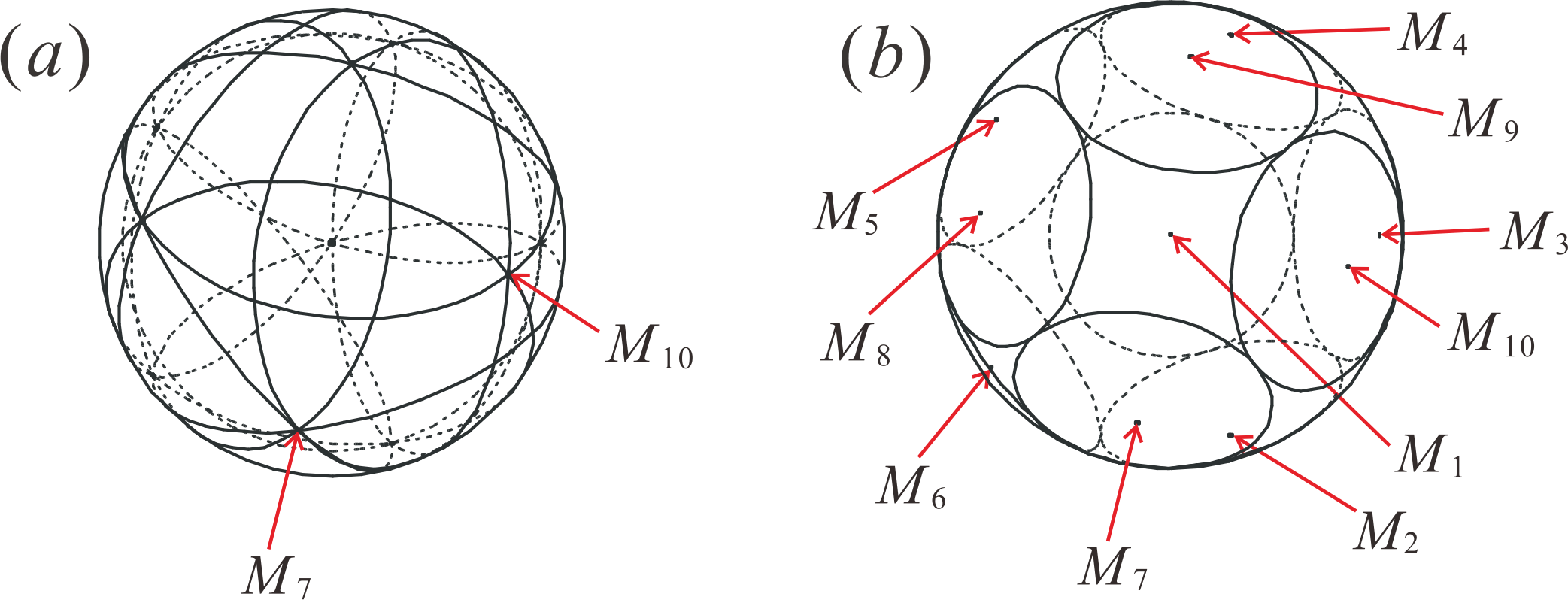}
\caption{{\protect\small ($a$) Our sequential covering for $N=10$. 
($b$) Our solution of Tammes problem for $N=10$. 
Both viewpoints are $(0, \;0, \;10)$. 
In this example, $M_{1}=(0, 0, -1)$, 
$M_{2} \approx (0.26335, -0.87585, -0.40439)$, 
$M_{3} \approx (0.91458, 0, -0.40439)$, 
$M_{4} \approx (0.26335, 0.87585, -0.40439)$, 
$M_{5} \approx (-0.76292, 0.50440, -0.40439)$, 
$M_{6} \approx (-0.77575, -0.57681, -0.25593)$, 
$M_{7} \approx (-0.13883, -0.78326, 0.60599)$, 
$M_{8} \approx (-0.79006, 0.092588, 0.60599)$, 
$M_{9} \approx (0.084546, 0.74290, 0.66405)$, 
and $M_{10} \approx (0.735778, -0.13295,$ 
$0.66405)$.
}
\label{fig_no1}
}
\end{figure}

We calculated the upper bound for $N = 2, \ldots ,12$ in our problem theoretically; 
the case $N=1$ is self-evident. As a result, we found the interesting fact 
that the solutions of our problem are strictly correspondent to those of 
the Tammes problem for $N = 2, \ldots ,12$~\cite{Sugi_Tane_2006, Sugi_Tane_2007}. 
Especially, as mentioned above, we succeeded to obtain the exact value for $N = 10$ 
(see (\ref{eq_r10}))~\cite{Sugi_Tane_2007}. 
In the case for $N=10$, when the centers are put on the spherical surface $\mathcal{S}$ 
according to our method, their arrangements are as shown in Figure~\ref{fig_no1}, 
for example. Hereafter, let us explain simply how to calculate the value of (\ref{eq_r10}).
When the centers of two spherical caps with angular radius $r$ are put 
respectively at $(0,\, 0,\, -1)$ and $\left( {\sin r, \,0, \,-\cos r}\right)$, 
we can obtain the coordinates $(x,y,z)$ of cross points where the perimeters 
of their spherical caps intersect by using simultaneous equations as follows:

\begin{equation}
\label{eq_2}
\left\{ {
\begin{array}{l}
-z=\cos r\,\text{,} \\ 
\sin r\cdot x-\cos r\cdot z=\cos r\,\text{,} \\ 
x^{2}+y^{2}+z^{2}=1\,\text{,}
\end{array}
}\right.
\end{equation}

\noindent
Solving (\ref{eq_2}), we obtain

\begin{equation}
\label{eq_3}
\left( {-\frac{\cos r\left( {\cos r-1}
\right) }{\sin r},\;\frac{\left( {\cos r-1}\right) \sqrt{2\cos r+1}}{\sin r}
,\;-\cos r}\right) \text{,}
\end{equation}

\begin{equation}
\label{eq_4}
\left( {-\frac{\cos r\left( {\cos r-1}
\right) }{\sin r},\;-\frac{\left( {\cos r-1}\right) \sqrt{2\cos r+1}}{\sin r}
,\;-\cos r}\right) \text{.}
\end{equation}

\noindent
In our method, taking into account Theorem1 in \cite{Sugi_Tane_2006}, the 
centers $M_1=(x_{m1}, \, y_{m1},\, z_{m1})$, $M_2=(x_{m2}, \, y_{m2},\, z_{m2})$, 
and $M_3=(x_{m3}, \, y_{m3},\,$ $z_{m3})$ are set at $(0,\; 0,\; -1)$, the coordinates 
of (\ref{eq_3}), and $( \sin r,\; 0,\; -\cos r )$, respectively. Here, let 
$\partial C_i$ be the perimeter of $C_{i} \; (\,i = 1, \ldots ,10\,)$. 
Next, according to our method, we choose $M_4=(x_{m4}, \, y_{m4},\, z_{m4})$ 
on the coordinates of (\ref{eq_4}). Then, $M_5=(x_{m5}, \, y_{m5},\, z_{m5})$ is put on 
one of the cross points of $\partial C_{4}$ and $\partial C_{1}$, and let it be 
outside $C_{3}$ (i.e., $( \sin r ( \cos ^{2}r-2\cos r-1)/ ( \cos r+1) ^{2},$ 
$2\cos r\,\sin r \sqrt{2\cos r+1}/ ( \cos r+1 ) ^{2},$ $-\cos r)$ ). In addition, 
we put $M_6=(x_{m6}, \, y_{m6},\,$ 
$z_{m6})$ on one of the cross points of 
$\partial C_{5}$ and $\partial C_{2}$, and let it be outside $C_{1}$ (i.e., 
$( 2\cos r\sin r ( \cos r-1)( 2\cos r+1) / (9\cos ^{3}r-\cos ^{2}r-\cos r+1),$  
$2\cos r \,\sin r ( \cos r-1) \sqrt{2\cos r+1}/ (9\cos ^{3}r-\cos ^{2}r-\cos r+1),$ 
$(-4\cos ^{4}r+\cos ^{3}r-5\cos ^{2}r-\cos r+1)/(9\cos^{3}r-\cos ^{2}r-\cos r+1 ) 
)$ ). 
At the time $C_{6}$ is put on a spherical surface $\mathcal{S}$, the 
uncovered region is reduced to a pentagon which is bounded by perimeters 
of spherical caps. 
Then, in order to put the centers of four spherical caps on the pentagonal 
uncovered region, we consider a spherical square of side-length $r$ on the 
pentagon. It is because, when the centers $M_7$, $M_8$, $M_9$, and $M_{10}$ are 
put on the vertices of spherical square of side-length $r$, the set 
$\cup _{\nu =1}^{10}C_{\nu }$ can cover the whole of $\mathcal{S}$ with our method.
Now, in our study, the centers $M_7$, $M_8$, $M_9$, and $M_{10}$ are put respectively 
on the cross points of $\partial C_{6}$ and $\partial C_{2}$, $\partial C_{7}$ 
and $\partial C_{6}$, $\partial C_{8}$ and $\partial C_{4}$, and $\partial C_{9}$ 
and $\partial C_{3}$. Note that the cross points chosen as arrangement of $M_7$, 
$M_8$, $M_9$, and $M_{10}$ are on the boundary of the pentagonal uncovered region.
Therefore, from the coordinates of $M_2$ and $M_6$, the coordinates of the 
centers $M_7 =(x_{m7}, \, y_{m7},\, z_{m7})$ are as obtained follows: 
 
\begin{eqnarray*}
&&\left( {\frac{\sin r\left( {\cos ^{3}r-5\cos ^{2}r-\cos r+1}\right) }{
9\cos ^{3}r-\cos ^{2}r-\cos r+1},\;-\frac{4\sin r\,\cos ^{2}r\;\sqrt{2\cos
r+1}}{9\cos ^{3}r-\cos ^{2}r-\cos r+1}}\right. ,\\ 
&&\left. {\quad -\frac{\cos r\left( {\cos ^{3}r+11\cos ^{2}r-\cos r-3}
\right) }{9\cos ^{3}r-\cos ^{2}r-\cos r+1}}\right) \text{.}
\end{eqnarray*}

\noindent
Next, we calculate the coordinates of $M_{10}$ without using the
coordinates of $M_{7}$. Here, let $s_{i \cdot j}$ denote the 
spherical distance between $M_{i}$ and $M_{j}$. Then, by applying 
the spherical cosine theorem to the spherical isosceles triangle 
$M_{6}M_{7}M_{10}$ of legs $s_{6 \cdot 7} =s_{7 \cdot 10} =r$, we have

\begin{equation*}
\cos (s_{6 \cdot 10}) 
=\frac{3\cos ^{3}r+2\cos ^{2}r-\cos r-2\left( {1-\cos ^{2}r}
\right) \sqrt{\cos r+2\cos ^{2}r}}{\left( {1+\cos r}\right) ^{2}}\text{.}
\end{equation*}

\noindent It is because the inner angle at $M_{6}$ of the spherical isosceles 
triangle $M_{6}M_{7}M_{10}$ is the sum of the interior angles of spherical 
equilateral triangle $M_{6}M_{7}M_{8}$ and spherical square $M_{7}M_{10}M_{9}M_{8}$.
In this connection, the inner angle of spherical equilateral triangle of side-length 
$r$ is $\cos ^{ - 1}( \cos r/(\cos r + 1) )$, 
and the inner angle of spherical square of side-length $r$ is 
$\cos ^{ - 1}( (\cos r- 1)/(\cos r + 1) )$.  
As a result, we can obtain the coordinates of $M_{10}=(x_{m10},
\, y_{m10},\, z_{m10})$ as a function of $r$ through the simultaneous 
equations as follows:

\begin{equation*}
\left\{ {
\begin{array}{l}
x_{m3}\cdot x_{m10}+y_{m3}\cdot y_{m10}+z_{m3}\cdot z_{m10}=\cos r\,\text{,} \\ 
x_{m6}\cdot x_{m10}+y_{m6}\cdot y_{m10}+z_{m6}\cdot z_{m10} =
\frac{3\cos ^{3}r+2\cos ^{2}r-\cos r-2\left( {1-\cos ^{2}r}\right) 
\sqrt{\cos r+2\cos ^{2}r}}{\left( {1+\cos r}\right) ^{2}} \text{,} \\ 
x_{m10} ^2 + y_{m10} ^2 +z_{m10} ^2=1\,\text{.}
\end{array}
}\right.
\end{equation*}

\noindent
Note that, in this report, the coordinates of $M_{10}$ are omitted since they 
are too complicated. Further, we get the equation of the following type 

\begin{equation}
\label{eq_5}
x_{m7}\cdot x_{m10}+y_{m7}\cdot y_{m10}+z_{m7}\cdot z_{m10}=\cos r.
\end{equation}

\noindent
When the equation (\ref{eq_5}) is solved against $r$ by using mathematical software, 
we obtained the value of (\ref{eq_r10}).

Danzer have solved the packing problem for $N=10$ through the consideration 
on irreducible graphs obtained by connecting those points, among $N$ points, 
whose spherical distance is exactly the minimal distance~\cite{Danzer_1963}. 
Then he needed the independent considerations for $N=10$. On the other hand, 
our systematic method is able to obtain a solution for $N$ by using the results 
for the case $N-1$ or $N-2$. In addition, we have considered the packing problem 
from the standpoint of sequential covering. The advantages of our approach 
are that we only need to observe uncovered region in the process of packing 
and that this uncovered region decreases step by step as the packing proceeds.

By using our method, the solutions of Tammes problem were obtained for 
$N = 2, \ldots ,12$~\cite{Sugi_Tane_2006, Sugi_Tane_2007}.
This fact is interesting and it is important that the exact value for 
$N= 10$ is found ~\cite{Sugi_Tane_2007}.

\bigskip
\noindent
{\textbf{Acknowledgments.}
The authors would like to thank Emeritus Professor L.~Danzer, University of Dortmund, 
and Emeritus Professor H.~Maehara, University of the Ryukyus, for their 
helpful comments.
The research was partly supported by the Grant-in-Aid for Scientific Research 
(the Grant-in-Aid for JSPS Fellows) from the Ministry of Education, Culture, 
Sports, Science, and Technology (MEXT) of Japan.


\begin{thebibliography}{9}


\bibitem{Danzer_1963}
L.~Danzer, \textit{Endliche punktmengen auf der 2-sph\"{a}re mit m\"{o}glichst 
gro{\ss}en Minimalabstand}. Universit\"{a}t G\"{o}ttingen, 1963. 
(Finite point-set on $S^{2}$ with minimum distance as large as possible,  
\textit{Discrete Mathematics}, \textbf{60} (1986) 3--66.)

\bibitem{Fejes_1972}
L.~Fejes~T\'{o}th, \textit{Lagerungen in der Ebene, auf der Kugel und im 
Raum}. Second edition. Springer-Verlag, Heidelberg, 1972. (Japanese trans. by I.~Higuchi, 
M.~Tanemura, \textit{Haichi no mondai --heimen$\cdot$kyuumen$\cdot$kuukan ni okeru}. 
Misuzusyobou, Tokyo, 1983).


\bibitem{Musin_Tarasov_2012}
O.R.~Musin, A.S.~Tarasov, The strong thirteen spheres problem,
\textit{Discrete \& Computational Geometry}, \textbf{48} (2012) 
128--141.
Available online: \url{http://dx.doi.org/10.1007/s00454-011-9392-2},
\url{http://arxiv.org/abs/1002.1439}
(accessed on 9 September 2015).


\bibitem{Musin_Tarasov_2015}
\textemdash, The Tammes problem for $N=14$, 
\textit{Experimental Mathematics}, \textbf{24} 
(2015) 460--468.
Available online: 
\url{http://www.tandfonline.com/doi/full/10.1080/10586458.2015.1022842}, 
\url{http://arxiv.org/abs/1410.2536}
(accessed on 9 September 2015).

\bibitem{Sugi_Tane_2006}
T.~Sugimoto, M.~Tanemura, Packing and Minkowski covering of congruent spherical 
caps on a sphere for $N = 2, \ldots,9$, \textit{Forma}, \textbf{21} (2006) 197--225.
Available online: \url{http://www.scipress.org/journals/forma/pdf/2103/21030197.pdf}
(accessed on 30 August 2015).

\bibitem{Sugi_Tane_2007}
T.~Sugimoto, M.~Tanemura, Packing and Minkowski covering of congruent spherical 
caps on a sphere, II: Cases of $N = 10, 11$ and $12$, \textit{Forma}, \textbf{22} 
(2007) 157--175.
Available online: \url{http://www.scipress.org/journals/forma/pdf/2202/22020157.pdf}
(accessed on 30 August 2015).

\bibitem{Tammes_1930}
P.M.L.~Tammes, On the origin of number and arrangement of the places of 
exit on the surface of pollen grains, \textit{Recueil de Travaux Botaniques 
N\'{e}erlandais}, \textbf{27} (1930) 1$-$84.

\bibitem{Teshi_Ogawa_2000}
Y.~Teshima, T.~Ogawa, Dense packing of equal circle on a sphere by
the minimum-zenith method: Symmetrical Arrangement, \textit{Forma}, 
\textbf{15} (2000) 347--364.
Available online: \url{http://www.scipress.org/journals/forma/pdf/1504/15040347.pdf}
(accessed on 30 August 2015).


\end{thebibliography}
\end{document}